\documentclass[A4paper, conference]{ieeeconf}
\IEEEoverridecommandlockouts
\usepackage{amsmath,amssymb,amsfonts,theorem}
\usepackage{graphicx,epsfig,color}
\usepackage{subfig}
\usepackage{algorithmic}
\usepackage{url}
\usepackage{bbm}

\graphicspath{{./figures/}}

\newcommand{\until}[1]{\{1,\dots, #1\}}

\renewcommand{\theenumi}{(\roman{enumi})}

\def\1{\mathbf{1}}

\renewcommand{\baselinestretch}{1}

\newcommand\oprocendsymbol{\hbox{$\square$}}
\newcommand\oprocend{\relax\ifmmode\else\unskip\hfill\fi\oprocendsymbol}

\def\Re{\mathcal{R}}
\def\P{\mathbb{P}}
\def\N{\mathbb{N}}
\def\R{\mathbb{R}}
\def\EE{\mathbb{E}}

\newcommand{\argmin}{\mathop{\operatorname{argmin}}}

\newtheorem{theorem}{Theorem}[section]
\newtheorem{proposition}[theorem]{Proposition}
\newtheorem{corollary}[theorem]{Corollary}
\newtheorem{definition}[theorem]{Definition}
\newtheorem{lemma}[theorem]{Lemma}
\newtheorem{assumption}[theorem]{Assumption}
\newtheorem{conjecture}[theorem]{Conjecture}
{\theorembodyfont{\rmfamily} 
\newtheorem{remark}[theorem]{Remark}

\newtheorem{example}[theorem]{Example}}

\def\log{\mathop{\operatorname{log}}}

\title{\bf An exponential-rate consensus-based
algorithms for estimation from relative measurements:
implementation and performance analysis
 \thanks{This work is supported by the European Community's Seventh Framework Programme [FP7/2007-2013] under
 grant agreement n. 257462 HYCON2 Network of excellence.}}

\author{ Ruggero Carli \and Andrea Carron \and Luca Schenato \and Marco Todescato%
  \thanks{R. Carli, A. Carron,  L. Schenato and M. Todescato are with the
    Department of Information Engineering, University of Padova, 
    Via Gradenigo 6/a, 35131 Padova, Italy 
    {\tt\small \{carlirug|carronan|schenato|todescat\}
    @dei.unipd.it}. 
  }}

\begin{document}
 \maketitle \thispagestyle{empty}
\pagestyle{empty}

\begin{abstract}

In this work we address the problem of optimal estimating the position of each agent in a network from relative noisy vectorial distances with its neighbors. Although the problem can be cast as a standard least-squares problem, the main challenge is to devise scalable algorithms that allow each agent to estimate its own position by means of only local communication and bounded complexity, independently of the network size and topology. We first propose a synchronous consensus-based algorithm that is guaranteed to have exponentially convergence rate to the optimal centralized least-squares solution. We then extend this algorithm to more realistic asynchronous implementations  via the use of local memory variables. We show that exponential convergence is still guaranteed under both uniform persistent communication and randomized persistent communication protocols. We finally complement the analytical results with some numerical simulations. 

\end{abstract}

\section{Introduction}

The proliferation of relatively inexpensive devices capable of communicating, computing, sensing, interacting with the environment and storing information  is promising an unprecedented number of novel applications throughout the cooperation of these devices toward a common goal. These applications include swarm robotics, wireless sensor networks, smart energy grids, smart traffic networks, smart camera networks. These applications also pose new challenges, of which \emph{scalability} is one of the major ones. Scalability is intended as the ability for an application to continue functioning without any dramatic performance degradation even if the number of devices involved keep increasing. In particular, an application is scalable if it is not necessary to increase HW resources or to adopt a more complex SW algorithms in each device even if the total number of devices increases.   

In this work we address the problem of designing algorithms that are capable to reconstruct the optimal estimate of the location of a device from noisy relative measurements from its neighbors in a connected network. In particular, we want to design distributed algorithms that allows each device to reconstruct its own position only from exchanging information with its neighbors, regardless of the size of the network. Moreover, these algorithms must be scalable, i.e. their computational complexity, bandwidth and memory requirements should be independent of the network size. These algorithms will be obtained as the solution of an optimization problem.

Distributed optimization has being attracting ever growing attention in the past years since many problems in large scale network have been cast as convex optimization problem. In particular, a large class of problem can be cast as the solution of the following optimization problem
\begin{equation}\label{opt1} \begin{array}{lll}x^*=&\mathrm{argmin}_{x_1,\ldots,x_N} & \sum_{i=1}^N f_i(x_i)\\ &s.t. & x_i=x_j, \forall i,j\end{array} \end{equation}
where $N$ is the number of nodes, $x_i\in\mathbb{R}^m$, and $f_i$ are convex functions. The function $f_i$ represents local cost of each agent $i$, but each agent must compute the minimizer $x^*$ of the sum all local costs. The local cost functions are separable but the additional contraint that all local variables $x_i$ must be the same, makes the optimization problem coupled. Many problems can be cast in this terms such as distributed least squares \cite{Xiao:05,Bolognani:10}, map building \cite{Schwager:09}, network utility maximization \cite{Johansson:06}, distributed learning and support vector machines  \cite{Boyd:11}. In these context several approaches have been proposed such as the distributed subgradient methods (SDMs) \cite{Nedich:09}, the alternating direction method of multipliers (ADMM)\cite{Boyd:11}, the Newton-Raphson consensus \cite{Zanella:11}, the control-based methods \cite{Wang:10}, and the Fast-Lipshitz \cite{Fischione:11}, and the distributed LP \cite{Notarstefano:11}.

The problem at hand in this work is of a different type and can be cast as the following unconstrained optimization problem: 
\begin{equation}\label{opt2} \min_{x_1,\ldots,x_N}  \sum_{(i,j)\in\mathcal{E}}^{|\mathcal{E}|} f_{ij}(x_i-x_j) \end{equation}
where $x_i\in \mathbb{R}^\ell$, $\mathcal{E}$ represents all the pair of nodes for which are available relative measurements and $f_{ij}$ are convex functions. Differently from the previous optimization problem, the local variables are unconstrained, but the cost functions are now not separable, and therefore the problem is once again coupled. Many problems can be cast in this framework such as sensor localization \cite{BarooahICISIP05,BarooahPhDThesis}, sensor calibration \cite{Bolognani:10}, clock synchronization \cite{Solis:06} and camera localization \cite{BorraACC12,TronCDC09}. 
 For example, in the context of localization from vectorial relative distance in a plane, the cost function $f_{ij}$ are given by:
$$f_{ij}(x_i-x_j)=||x_i-x_j-z_{ij}||^2$$
where $z_{ij}\in\mathbb{R}^\ell$ is the noisy measurement of the relative (vector) distance of node $i$ from node $j$. As a consequence, the optimization problem in Eqn.~\eqref{opt2} becomes a distributed least-square problem that in principle could be cast as the optimization problem in Eqn.~\eqref{opt2}. However, in this case, each node will need to compute the location of all other nodes, i.e. the size of the local variable becomes of size $m=N\ell$, which according to our objective is not scalable. There exist alternative approaches that try to exploit the special structure of the problem, but they either show oscillatory behavior as in consensus-based with constant weights \cite{Bolognani:10}, or the convergence rate decreases only as $1/k$ as in the Randomized Kaczmarz with Under-Relaxation \cite{Freris:12}. This slow convergence rate is mainly due to the fact that these algorithms adopt asynchronous gradient-based algorithms whose step-size decreases to zero as progress.

The contribution of this work is to provide both synchronous and asynchronous algorithms which are scalable and have proven exponential rate of convergence under mild assumptions. The main idea is to cast the problem as a consensus problem under some suitable change of coordinates, and then solve the problem by applying synchronous broadcast consensus iterations. In a more realistic asynchronous broadcast communication, the problem is solved similarly by adding some extra memory variables at each node which keep track of the estimated location of their neighbors, i.e. the nodes from which they collected the relative distance measurements. Estimates of these local variables eventually converge to the estimates of the neighbors, thus guaranteeing the convergence of the whole algorithm, at the price of some delay.    


\section{Problem Formulation}
The problem we deal with is that of estimating $N$ variables $x_1,\ldots, x_N$ from noisy measurements of the form 
\begin{equation}\label{eq:MeasurementModel}
z_{ij}:=x_i-x_j+n_{ij}, \qquad i,j \in \until{N},
\end{equation}
where $n_{ij}$ is zero-mean measurement noise. Though the variables
are often vector-valued, for simplicity, in this paper we assume that $x_i \in \R$, $i \in \until{N}$. 

This estimation problem can be naturally associated with a
\emph{measurement graph} $\mathcal{G} = (V ; \mathcal{E})$. The vertex set $V$ of the measurement graph
consists of the set of nodes $V = \left\{1,\ldots,N\right\}$ where $N$ is the number of nodes,
while its edge set $\mathcal{E}$ consists of all of the ordered pairs of nodes $(i,j)$ such that
a noisy measurement of the form \eqref{eq:MeasurementModel}
between $i$ and $j$ is available to node $i$. The measurement
errors on distinct edges are assumed uncorrelated. 
The measurement graph $\mathcal{G}$ is a directed
graph since $(i, j) \in  \mathcal{E}$ implies the measurement $z_{ij}$ is available to node $i$, 
while $(j, i) \in  \mathcal{E}$  implies the measurement $z_{ji}$ is available to node $j$, and these two are in general distinct.

To formally state the problem we aim at solving we need some preliminary definitions. 

Assume $|\mathcal{E}| =M$ and let us introduce the \emph{incidence matrix} $A \in \R^{M \times N}$ of $\mathcal{G}$ defined as $A=[a_{ei}]$, where $a_{ei}=1, -1, 0$, if edge $e$ is incident on node $i$ and directed away from i, is incident on node $i$ and directed toward it, or is not incident on node $i$, respectively.\\

Now, let $\textbf{x} \in \R^N$ be the vector obtained stacking together all the variables $ x_1, \ldots, x_N$, i.e., ${\bf x}=[x_1,\ldots, x_N]^T$, where given a vector $v$ with $v^T$ we denote its transpose, and let $\textbf{z} \in \R^{M}$ and $\textbf{n} \in \R^M$ be the vectors obtained stacking together all the measurements $z_{ij}$ and the noises $n_{ij}$, respectively. 
Additionally, let $R_{ij}>0$ denote the covariance of the zero mean error $n_{ij}$, i..e, $R_{ij}=\EE[n^2_{ij}]$,  where $\EE$ denotes the expectation operator, and let $R\in \R^{M \times M}$ be the diagonal matrix collecting in its diagonal the covariances of the noises $n_{ij}$, $(i,j)\in \mathcal{E}$, i.e., $R=\EE[\bf{n}\bf{n}^T]$. Finally let $\1$ be the column vector with all components equal to one. 

Observe that equation $\eqref{eq:MeasurementModel}$ can be rewritten in a vector form as
$$
{\bf z} = A  {\bf x}+ \bf n
$$ 
Define the set
\begin{equation*}
\chi := \argmin_{{\bf  x} \in \R^N} \,  ({\bf z} - A {\bf x})^T  R^{-1} ({\bf z} - A \bf x).
\end{equation*}
The goal is to construct an optimal estimate $x^*$ of $\bf x$ in a least square sense, namely, to compute
\begin{equation}\label{eq:least-square}
x^* \in \chi   
\end{equation}
Assume the measurement graph $\mathcal{G}$ to be \emph{weakly connected}\footnote{Given a directed graph, a path from a node to another node that does not respect the orientation of the edges is called an \emph{undirected path}. A directed graph is said to be \emph{weakly connected} if there is a undirected path form any node to any other node.}, then
it is well known (see \cite{BarooahPhDThesis}) that
$$
\chi = \left\lbrace \left(A^T R^{-1}A\right)^{\dagger}A^T R^{-1} \mathbf{z} + \alpha\1 \right\rbrace.
$$
Moreover let
$$
x_{\text{opt}}^* = \left(A^T R^{-1}A\right)^{\dagger}A^T R^{-1} \mathbf{z},
$$
then $x_{\text{opt}}^*$ is the minimum norm solution of \eqref{eq:least-square}, i.e., 
$$
x_{\text{opt}}^*=\min_{x^* \in \chi} \, \parallel x^* \parallel
$$
The matrix $A^T R^{-1}A$ is called in literature the \emph{Weighted Generalized Grounded Laplacian} (\cite{BarooahPhDThesis}).

\begin{remark}
Observe that, just with relative measurements, determining the $x'_i$s
is 
only possible up to an additive constant. This
ambiguity might be avoided by assuming that a node (say node $1$) is used as
reference node, i.e., $x_1 = 0$.
\end{remark}

\section{A synchronous distributed consensus based solution}
\label{sec:SynchImplementation}

To compute an optimal estimate $x^*$ directly, one needs all the measurements and their covariances $(\textbf{z}, R)$, and the topology of the measurement graph $\mathcal{G}$. 
In this section the goal is to compute the optimal solution in a distributed fashion, employing only local communication. In particular we assume that a node $i$ and another node $j$ can communicate with each other if either $(i,j) \in \mathcal{E}$ or $(j,i) \in \mathcal{E}$. Accordingly a node $i$ is said to be a neighbor of another node $j$ (and vicecersa) if either $(i,j) \in \mathcal{E}$ or $(j,i) \in \mathcal{E}$. For $i \in \until{N}$, by $\mathcal{N}_i$ we denote the set of neighbors of node $i$, namely
$$
\mathcal{N}_i=\left\{j \in V\,\,  \text{such that either $(i,j) \in \mathcal{E}$ or $(j,i) \in \mathcal{E}$}\right\}.
$$

In what follows we introduce a distributed solution which is based on standard linear consensus algorithm. A discussion of the linear consensus algorithm can be found in the review paper \cite{Olfati:07}, hence we refrain from describing it here. Instead we make the presentation of the algorithm self-contained. At first we assume that the communications among the nodes are synchronous, namely the nodes perform their transmissions' actions at the same instant times, and design the algorithm for that scenario. We refer to this algorithm as the \emph{synchronous consensus-based localization} algorithm (denoted hereafter as s-CL algorithm).   
In section \ref{sec:AsynchImp} we modify the s-CL algorithm to make it suitable to asynchronous communications. 

We assume that before running the s-CL algorithm, the nodes exchange with their neighbors their relative measurements as well as the associated covariances. So every node has access to the measurements on the edges that are incident to it, whether the edge is directed to or away from it. Each node uses the measurements obtained initially for all future computations. 

The s-CL algorithm is formally described as follows. 
\begin{LaTeXdescription}
\item[Processor states:] For $i\in \until{N}$, node $i$ stores in memory the measurements $\left\{z_{ij}, (i,j) \in \mathcal{E}\right\}$, $\left\{z_{ji}, (j,i) \in \mathcal{E}\right\}$, and the associated covariances $\left\{R_{ij}, (i,j) \in \mathcal{E}\right\}$, $\left\{R_{ji}, (j,i) \in \mathcal{E}\right\}$. Moreover node $i$ stores in memory also an estimate $\hat{x}_i$ of $x_i$. 
\item[Initialization:] Every nodes initializes its estimate to an arbitrary value.
\item[Transmission iteration:] For $k\in \N$, at the start of the $(k+1)$-th iteration of the algorithm, node $i$ transmits its estimate $\hat{x}_i(k)$ to all its neighbors. It also gathers the $k$-th estimates of its neighbors, $\hat{x}_j(k)$, $j \in \mathcal{N}_i$.  
\item[Update iteration:] For $k\in \N$, node $i$, $i \in \left\{1,\ldots,N\right\}$, based on the information received from its neighbors, updates its estimate as follows
$$
\hat{x}_i(k+1):= p_{ii} \,\hat{x}_i(k) + \sum_{j \in \mathcal{N}_i}  p_{ij} \, \hat{x}_j(k)+ b_i
$$
where 
$$
b_i=\epsilon \sum_{(i,j)\in \mathcal{E}} R^{-1}_{ij} z_{ij} - \epsilon \sum_{(j,i)\in \mathcal{E}} R^{-1}_{ji} z_{ji}
$$
and where
$$
p_{ij}=\left\{
\begin{array}{ccc}
 \epsilon (R^{-1}_{ij}+R^{-1}_{ji}) & \text{if $(i,j) \in \mathcal{E}$ and $(j,i) \in \mathcal{E}$} \\
\epsilon  R^{-1}_{ij} & \text{if $(i,j) \in \mathcal{E}$ and $(j,i) \notin \mathcal{E}$}\\
\epsilon R^{-1}_{ji} & \text{if $(j,i) \in \mathcal{E}$ and $(i,j) \notin \mathcal{E}$}\\
\end{array}
\right.
$$
and
$$
p_{ii}=1-\sum_{j\in \mathcal{N}_i} p_{ij}
$$
being $\epsilon$ a positive constant \emph{a-priori} assigned to the nodes.
\end{LaTeXdescription}
Now, let $P \in R^{N \times N}$ be the matrix defined by the weights $p_{ij}$ above introduced. One can see that such matrix $P$ is equal to 
$$
P = I-\epsilon A^T R^{-1} A
$$ 
Moreover let
$$
b=
\epsilon A^TR^{-1} \textbf{z}
$$
Then the s-CL algorithm can be written in a compact form as
$$
\hat{x}(k+1)=P \hat{x}(k) + b
$$
To characterize the convergence properties of the s-CL algorithm, we introduce some definitions. Firstly, for $i \in \until{N}$, let $d_{max}= \max \left\{|\mathcal{N}_i|, i\in \until{N}\right\}$. Secondly let $R_{min} = \min \left\{R_{ij}, (i,j) \in \mathcal{E}\right\}$. 
Observe that if $0< \epsilon < 1/(2 d_{\text{max}} R^{-1}_{min})$ then the matrix $P$ is stochastic. Moreover if the graph $\mathcal{G}$ is weakly connected then the matrix $P$ is primitive, i.e., it  satisfies the following two properties
\begin{enumerate}
\item $P$ has only one eigenvalue equal to $1$, $P\1=\1$ and $\1^TP= \1^T$, since $P$ is double stochastic;
\item all the other eigenvalues of $P$ are all strictly inside the unitary circle. 
\end{enumerate}
We recall that, based on standard results on nonnegative matrices and linear consensus algorithms, the above two properties imply $\lim_{k \to \infty} P^k =\frac{\1\1^T}{N}$.  Moreover let $\rho_{ess}$ denote the essential spectral radius of $P$, namely, the second largest in absolute value eigenvalue of $P$ (see \cite{RC-FF-AS-SZ:08}). Clearly $0 \leq \rho_{ess} < 1$.
We have the following Proposition.

\begin{proposition}\label{prop:1}
Let $\epsilon$ be such that $0< \epsilon < 1/(2 d_{\text{max}} R^{-1}_{min})$. Moreover let $\hat{x}_i$, $i \in \left\{1,\ldots,N\right\}$, be initialized to any real number. Let $\mathcal{G}$ be \emph{weakly connected}. Consider the s-CL algorithm running over $\mathcal{G}$. 
Then the following two facts hold true
\begin{enumerate}
\item the evolution $k \to \hat{x}(k)$ asymptotically converges to an optimal estimate $x^*\in \chi$, i.e., there exists $\alpha \in \R$, such that
$$
\lim_{k \to \infty} \hat{x}(k)=x_{\text{opt}}^* + \alpha\1;
$$ 
where $\alpha$ linearly depends on $x(0)$.
\item the convergence is exponential, namely, there exists $C>0$  such that
$$
\| \hat{x}(k)-\left(x_{\text{opt}}^* +\alpha\1\right)\| \leq C \rho_{ess}^k \| \hat{x}(0)-\left(x_{\text{opt}}^* +\alpha\1\right)\|.
$$ 
\end{enumerate}
\end{proposition}

\begin{proof}
We start by proving item (i).
Let us define the change of variable $\xi = \hat{x}-x_{\text{opt}}^*$. Since 
$$
x_{\text{opt}}^* = Px_{\text{opt}}^* + b
$$
it is possible to write
\begin{align*}
\hat{x}(k+1)-x_{\text{opt}}^* &= P\hat{x}(k) + b - x_{\text{opt}}^* \\
&= P\hat{x}(k) + b - (Px_{\text{opt}}^* + b) \\
&=P(\hat{x}(k)-x_{\text{opt}}^*)
\end{align*}
and, in turn, 
$$
\xi(k+1) = P\xi(k).
$$
The above equation describes the iteration of the classical consensus algorithm.
Since $P$ is a primitive doubly stochastic matric, we have that
$$
\xi(k) \rightarrow \frac{\1\1^T}{N}\xi(0)
$$
where $\xi(0)=x(0)-x_{\text{opt}}^*$. This implies that
$$
x(k) \rightarrow x_{\text{opt}}^* + \frac{\1\1^T}{N}x(0) - \frac{\1\1^T}{N}x_{\text{opt}}^*
$$
Being $\frac{\1\1^T}{N}x^*_{\text{opt}}=0$ this prove the result.

Concerning item (ii) it is well known (\cite{RC-FF-AS-SZ:08}) that the convergence rate of a consensus algorithm ruled by a primitive matrix $P$, is exponential and is upper bounded by the essetial spectral radius $\rho_{ess}$.
\end{proof}

\begin{remark}
The s-CL algorithm is similar to the algorithm proposed in \cite{FP:FF:RW:12}. However in \cite{FP:FF:RW:12}, the measurement graph is assumed to be undirected, namely, both measurements $z_{ij}$ and $z_{ji}$ are available to node $i$ and $j$ under the additional assumption that $z_{ij}=-z_{ji}$. 
\end{remark}
\begin{remark}
The authors in \cite{Barooah:07b} solved the problem formulated in \eqref{eq:least-square} proposing a synchronous algorithm that implements the Jacobi iterative method. The performance of this algorithm, in terms of rate of convergence to the optimal solution, is similar, for many families of measurement graphs, to the one of the synchronous consensus-based algorithm introduced in this section. However, to the best of our knowledge, there is no an asynchronous implementation of the Jacobi-based algorithm proposed in  \cite{Barooah:07b}, which has been shown to be provably convergent to the optimal solution. 
\end{remark}


\section{An asynchronous implementation of distributed consensus based solution}
\label{sec:AsynchImp}

The distributed algorithm illustrated in the
previous section, has an important limitation: it is applicable
only to sensors' networks with synchronized and reliable communication.
Indeed, the s-CL algorithm requires that there
exists a predetermined common communication schedule
for all nodes and, at each communication round, each
node must simultaneously and reliably communicates its
information.
The aim of this section is to reduce the communication
requirements of the s-CL algorithm, in particular in terms of synchronization.
To do so, we next introduce the \emph {asynchronous Consensus-based Localization} algorithm (denoted as a-CL hereafter).
This algorithm is based on an asymmetric broadcast
communication protocol. Differently from the s-CL, at each
iteration of the a-CL there is only one node transmitting
information to all its neighbors. 

Since the actual value of neighboring estimates are not available at each iteration, we assume that each node stores in its local memory a copy of the neighbors' variables recorded from the last communication they performed. 
 For $j \in \mathcal{N}_i$, we denote by $\hat{x}_j^{(i)}(k)$ the estimate of $x_j$ kept in $i$'s local memory at the end of the $k$-th iteration.  If node $j$ performed its last transmission to node $i$ during $h$-th iteration, $h \leq k$, then $\hat{x}_j^{(i)}(k)=\hat{x}_j(h)$.

The a-CL algorithm is formally described as follows. 

\begin{LaTeXdescription}
\item[Processor states:] For $i\in \until{N}$, node $i$ stores in memory the measurements $z_{ij}$, $z_{ji}$ and the covariances $R_{ij}$, $R_{ji}$ for all $j \in \mathcal{N}_i$. Moreover node $i$ stores in memory also the estimate $\hat{x}_i$ of $x_i$ and, for $j\in \mathcal{N}_i$ an estimate $\hat{x}_j^{(i)}$ of $\hat{x}_j$,
\item[Initialization:] Every node $i$ initializes its estimate $\hat{x}_i$ and the variables $\hat{x}_j^{(i)}$, $j\in \mathcal{N}_i$, to arbitrary values.
\item[Transmission iteration:] For $k\in \N$, at the start of the $(k+1)$-th iteration of the algorithm, there is only one node, say $i$, which
transmits information to its neighbors; precisely, node $i$ sends the value of its estimate $\hat{x}_i(k)$ to node $j$, $j \in \mathcal{N}_i$. 
\item[Update iteration:] 
For $j \in \mathcal{N}_i$, node $j$ performs the following actions in order
\begin{enumerate}
\item it sets $\hat{x}_i^{(j)}(k+1)=\hat{x}_i(k)$, while for $s\in \mathcal{N}_j \setminus \left\{i\right\}$, $\hat{x}_s^{(j)}$ is left unchanged, i.e., $\hat{x}_s^{(j)}(k+1)=\hat{x}_s^{(j)}(k)$;
\item it updates $\hat{x}_j$ as 
\begin{equation}\label{eq:UpdateAsyn}
\hat{x}_j(k+1):= p_{jj} \hat{x}_j(k) +\sum_{h \in \mathcal{N}_j} p_{jh} \hat{x}_h^{(j)}(k+1)+ b_j.
\end{equation}
\end{enumerate}
Clearly for $s \notin \mathcal{N}_i$, $\hat{x}_s$ is left unchanged during the $(k+1)$-th iteration of the algorithm, i.e, $\hat{x}_s(k+1)=\hat{x}_s(k)$.
\end{LaTeXdescription}

We characterize now the convergence properties of the
a-CL. We consider two different scenarios which we introduce in the following definition.

\begin{definition}
A network of $N$ nodes is said to be a \emph{uniformly persistent communicating network} if there exists a positive integer
number $\tau$ such that, for all $k \in \N$, each node transmits the value of its estimate to its neighbors at least once within the time interval $[k, k+\tau)$.  

A network of $N$ nodes is said to be a \emph{randomly persistent communicating network} if there exists a $N$-upla of probabilities $(\beta_1,\ldots, \beta_N)$ such that $\beta_i>0$, for all $i \in \until{N}$, and $\sum_{i=1}^N \beta_i=1$, and such that, for all $k \in \N$,
$$
\P \left[\text{the transmitting node at iteration $k$ is node $i$} \right]=\beta_i.
$$
\end{definition}

\subsection{Convergence of the a-CL algorithm under uniform persistent communications}

The following result characterizes the convergence properties of the a-CL when the network is uniformly persistent communicating. 

\begin{proposition}
Consider a \emph{uniformly persistent communicating network} of $N$ nodes with associated a {weakly connected measurement graph} $\mathcal{G}$ running the a-CL algorithm. 
Let $\epsilon$ be such that $0< \epsilon < 1/( 2 d_{\text{max}} R_{min}^{-1})$. Moreover let $\hat{x}_i$, $i \in \left\{1,\ldots,N\right\}$, $\hat{x}_j^{(i)}$, $j \in \mathcal{N}_i$, be initialized to any real number. Then the following facts hold true
\begin{enumerate}
\item the evolution $k \to \hat{x}(k)$ asymptotically converges to an optimal estimate $x^* \in \chi$, i.e., there exists $\alpha \in \R$ such that  
$$
\lim_{k \to \infty} \hat{x}(k)=x_{\text{opt}}^* +\alpha\1;
$$ 
\item the convergence is exponential, namely, there exists $C>0$ and $0 \leq \rho <1$ such that
$$
\| \hat{x}(k)-\left(x_{\text{opt}}^* +\alpha\1\right)\| \leq C \rho^k \| \hat{x}(0)-\left(x_{\text{opt}}^* +\alpha\1\right)\|.
$$ 
\end{enumerate}
\end{proposition}

\begin{proof}
Observe that, for $i \in \until{N}$ and $j \in \mathcal{N}_i$, we have $\hat{x}_i^{(j)}(k)=\hat{x}_i(\tau_i(k))$, where $\tau_i(k)<k$ denotes the iteration during which node $i$ has perfomed its last transmission up to iteration $k$. Hence Equation \eqref{eq:UpdateAsyn} can be rewritten as
\begin{equation}\label{eq:UpdateAsynDelay}
\hat{x}_j(k+1):= p_{jj} \hat{x}_j(k) +\sum_{h \in \mathcal{N}_j} p_{jh} \hat{x}_h(\tau_h(k)) + b_j.
\end{equation}
Observe that
$$
x_{\text{opt}}^* = Px_{\text{opt}}^* + b
$$ 
and, in particular, 
$$
\left[x_{\text{opt}}^*\right]_j= \left[Px_{\text{opt}}^*\right]_j + b_j = p_{jj}\left[x_{\text{opt}}^*\right]_j + \sum_{h\in\mathcal{N}_j} p_{jh}\left[x_{\text{opt}}^*\right]_h + b_j.
$$
By subtracting $\left[x_{\text{opt}}^*\right]_j$ to both the left-hand side and the right-hand side of \eqref{eq:UpdateAsynDelay} one gets
\begin{align*}
& \hat{x}_j(k+1) - \left[x_{\text{opt}}^*\right]_j \\
& \qquad = p_{jj}\hat{x}_j(k) + \sum_{h\in\mathcal{N}_j} p_{jh}\hat{x}_h(\tau_h(k)) + \\
&\qquad\qquad + b_j- \left( p_{jj}  \left[x_{\text{opt}}^*\right]_j + \sum_{h\in\mathcal{N}_j} p_{jh} \left[x_{\text{opt}}^*\right]_h + b_j \right)\\
&\qquad= p_{jj}(\hat{x}_j(k) -  \left[x_{\text{opt}}^*\right]_j) + \\
&\qquad\qquad +\sum_{h\in\mathcal{N}_j} p_{jh}\left(\hat{x}_h(\tau_h(k)) -  \left[x_{\text{opt}}^*\right]_h\right).
\end{align*}
Let us introduce the auxiliary variable $\xi(k) = \hat{x}(k) - x_{\text{opt}}^*$. From the above equation we can write
\begin{equation}
\label{eq:changeOfVar}
\xi_j(k+1) = p_{jj}\xi_j(k) + \sum_{h\in\mathcal{N}_j} p_{jh}\xi_h(\tau_h(k))
\end{equation}
Observe that \eqref{eq:changeOfVar} describes the evolution of a consensus algorithm in presence of uniformly bounded delays. Hence, by invoking Proposition $1$ in \cite{Nedic:08}, we can conclude that $\xi(k)$ converges exponentially to consensus, and, in turn, that $\hat{x}(k)$ converges exponentially to an optimal solution $x^*$.
\end{proof}

\subsection{Convergence of the a-CL algorithm under randomly persistent communications}

The following result characterizes the convergence properties of the a-CL when the network is randomly persistent communicating. 

\begin{proposition}
Consider a \emph{randomly persistent communicating network} of $N$ nodes with associated a {weakly connected measurement graph} $\mathcal{G}$ running the a-CL algorithm. 
Let $\epsilon$ be such that $0< \epsilon < 1/( 2 d_{\text{max}} R_{min}^{-1})$. Moreover let $\hat{x}_i$, $i \in \left\{1,\ldots,N\right\}$, $\hat{x}_j^{(i)}$, $j \in \mathcal{N}_i$, be initialized to any real number. Then the following facts hold true
\begin{enumerate}
\item
the evolution $k \to \hat{x}(k)$ converges almost surely to an optimal solution $x^*\in\chi$, i.e., there exists $\alpha \in \R$ such that
$$
\P\left[\lim_{k \to \infty} \hat{x}(k)=x^*_{\text{opt}} + \alpha\1\right]=1.
$$
\item the evolution $k \to \hat{x}(k)$ is exponentially convengent in mean-square sense, i.e., there exist $C>0$ and $0 \leq \rho <1$ such that
\begin{align*}
& \lim_{k \to \infty} \EE\left[\|\hat{x}(k)-(x^*_{\text{opt}} + \alpha\1)\|^2\right] \\
&\qquad \qquad\leq C  \rho^k \EE\left[\|\hat{x}(0)-(x^*_{\text{opt}} + \alpha\1)\|^2\right].
\end{align*}
\end{enumerate}
\end{proposition}

\begin{proof}
We start by rewriting the updating step of the a-CL in a more convenient way. Observe preliminarily that, under the assumption of reliable communications over the network, the broadcast protocol lets only two information about the estimate of $x_{i}$, $\forall i \in V$, to flow through the network. Specifically, $\hat{x}_{i}(k)$, that is the current value of the estimate $\hat{x}_{i}$ at iteration $k$, and $\hat{x}_{i}(t'_i(k))$, being $t'_i(k)$ the iteration during which node $i$ has performed its last transmission up to iteration $k$ of a-CL (that is, the value of $\hat{x}_i^{(i)}$ at its last communication round). Indeed notice that, for $j\in \mathcal{N}_i$, $\hat{x}_i^{(j)}(t'')= \hat{x}_i(t'_i(k))$ for all $t''$ such that $t'_i(k)< t'' \leq k$. Let us define $x'_i(k)=\hat{x}_i(k)$ and $x''_i(k)=\hat{x}_i(t'_i(k))$ and, accordingly, let $x'(k)=\left[x'_1(k),\ldots, x'_N(k)\right]^T$ and $x''(k)=\left[x''_1(k),\ldots, x''_N(k)\right]^T$. 
Morever let $Q_i \in \R^{2N \times 2N}$ be defined as
\begin{equation}\label{eq:Qi}
Q_i =\left[
\begin{array}{cc}
Q_{11}^{(i)} & Q_{12}^{(i)}\\
Q_{21}^{(i)} & Q_{22}^{(i)}
\end{array}
\right]
\end{equation}
where
\begin{eqnarray}
Q_{11}^{(i)} &=& \sum_{h\notin\mathcal{N}_{i}}e_{h}e_{h}^{T} + \sum_{j\in\mathcal{N}_{i}}\left(p_{jj}e_{j}e_{j}^{T} + p_{ji}e_{j}e_{i}^{T}\right)\nonumber\\
Q_{12}^{(i)} &=& \sum_{j\in\mathcal{N}_{i}}e_{j}\left(\sum_{h\in\mathcal{N}_{j}/i}p_{jh}e_{h}^{T}\right)\nonumber\\
Q_{21}^{(i)} &=& e_{i}e_{i}^{T}\nonumber\\ 
 Q_{22}^{(i)} &=& I - e_{i}e_{i}^{T} \nonumber
\end{eqnarray}
Observe that, for $i \in \until{N}$, $Q_i$ is a $2N$-dimensional stochastic matrix.   
Finally let 
$$
B_i=
\left[
\begin{array}{c}
\sum_{j \in \mathcal{N}_i} e_j^T b\\
0_N
\end{array} 
\right]
$$

Assume, without loss of generality, that node $i$ is the node performing the transmission during the $k+1$-th iteration of the a-CL. Hence the udpating step of the a-CL can be written in compact form as
\begin{equation}
\left[
\begin{array}{c}
x'(k+1)\\
x''(k+1)
\end{array}
\right]= Q_i 
\left[
\begin{array}{c}
x'(k)\\
x''(k)
\end{array}
\right] + B_i.
\end{equation}
Now let us introduce the auxiliary variable
$$
\xi(k) = 
\left[
\begin{array}{c}
x'(k)\\
x''(k)
\end{array}
\right]-
\left[
\begin{array}{c}
x_{opt}^*\\
x_{opt}^*
\end{array}
\right]
$$
By exploiting the fact that, for $i \in \until{N}$,
\begin{equation}
\left[
\begin{array}{c}
x_{opt}^*\\
x_{opt}^*
\end{array}
\right]= Q_i 
\left[
\begin{array}{c}
x_{opt}^*\\
x_{opt}^*
\end{array}
\right] + B_i
\end{equation}
we have that the variable $\xi$ satisfies the following $2N$-dimensional recursive equation
\begin{equation}\label{eq:EqXi}
\xi(t+1)=Q_i\xi(t).
\end{equation}
Observe that $\hat{x}(k) \to x_{\text{opt}}^*+\alpha \1$ if and only if $\xi(k) \to \alpha \1$. Moreover,
since $Q_i$ is a stochastic matrix for any $i\in \until{N}$, we have that \eqref{eq:EqXi} describes the evolution of a randomized consensus algorithm whose convergence properties can be analyzed following the treatment in \cite{FF-SZ:08a}.

Now let $\sigma$ be the random process such that $\sigma(k)$ denotes the node performing the transmission action at the beginning of the $k+1$-th iteration. Clearly, in the randomized scenario we are considering, we have that, for $i\in \until{N}$, $\P[\sigma(k)=i]=\beta_i$ for all $k$. Let
$$
S(k)=\prod_{h=0}^k Q_{\sigma(h)}.
$$
Observe that $S(k)$ inherits the same block structure of the matrices $\left\{Q_i\right\}_{i=1}^N$, namely we can write
$$
S(k)=
\left[
\begin{array}{cc}
S_{11}(k) & S_{12}(k)\\
S_{21}(k) & S_{22}(k)
\end{array}
\right]
$$

As consequence of Theorem $3.1$ in \cite{FF-SZ:08a} the a-CL reaches almost surely consensus if and only if, for every $i$ and $j$ in $V$
\begin{equation}\label{eq:ConditionRandomizedConsensus}
\P\left[\mathcal{E}_{ij}\right]=1,
\end{equation}
where 
$$
\mathcal{E}_{ij}=\left\{\exists \ell, \exists k \,|\, S_{i\ell}(k)S_{j \ell}(k)>0 \right\}.
$$
Now observe that, since the measurement graph is weakly connected, then the communication graph is a connected undirected graph. This fact together with the fact the diagonal elements of $Q^{(i)}_{11}$ are all positive for any $i\in \until{N}$ implies that there exists almost surely $\bar{k}$ such that, for all $k'\geq \bar{k}$, all the elements of the matrix $S_{11}(k')$ are strictly greater than $0$. Assume now, without loss of generality, that $\sigma(k')=i$, for $k'\geq k$. Then, since the $i$-th row of $S_{21}(k'+1)$ is equal to $e_i e_i^T S_{11}(k')$, it turns out that, all the elements of the $i$-th row of $S_{21}(k'+1)$ are strictly greater than $0$. Moreover, it is easy to see that they will remain stricty greater than $0$ also for any $k''\geq k'$. Hence we can argue that, there exists almost surely, also a $\bar{k}'$ such that  for all $k'\geq \bar{k}'$, all the elements of the matrix $S_{21}(k')$ are strictly greater than $0$. It follows that the property stated in \eqref{eq:ConditionRandomizedConsensus} is satisfied for any $k\geq \bar{k}'$ and for any $\ell \in \until{N}$. This concludes the proof of item (i).

Concerning item (ii), we again resort to the results in \cite{FF-SZ:08a}. Let $\Omega=I-1/2N \1\1^T$ where in this expression we assume that $I$ is the $2N$-dimensional identity matrix and the vector $\1$ is $2N$-dimensional. From the results in \cite{FF-SZ:08a}, it follows that to study the rate of convergence of $\EE\left[\|\xi(k)- \alpha\1 \|^2\right]$ is equivalent to study the convergence rate of $\EE \|\Omega \xi(k)\|^2$ and in particular of the linear recursive system 
$$
\Delta(t+1)= \EE \left[Q_{\sigma(0)}^T \Delta(t)Q_{\sigma(0)} \right]
$$
where $\Delta(0)=\Omega$. Observe that $\Delta(t)$ is the evolution of a linear dynamical system which can be written in the form 
$$
\Delta(t+1)=\mathcal{L}(\Delta(t))
$$
where $\mathcal{L}: \R^{2N \times 2N} \to \R^{2N \times 2N}$ is given by
$$
\mathcal{L}(M)=\EE\left[Q_{\sigma(0)}^T M Q_{\sigma(0)} \right].
$$
As highlighted in \cite{FF-SZ:08a}, the linear operator $\mathcal{L}$ can be represented by the matrix $\bf{L}=\EE[Q_{\sigma(0)} \otimes Q_{\sigma(0)}]^T$ where $\otimes$ denotes the Kronecker product of matrices. Following the proof of Proposition $4.3$ of  \cite{FF-SZ:08a}, one can see that $\bf{L}^T$ is a primitive stochastic matrix which, therefore, has the eigenvalue $1$ with algebraic multiplicity $1$. Moreover, $\bf{L}^T(\1 \otimes \1)=(\1 \otimes \1)$ and $(\1 \otimes \1) (\Omega \otimes \Omega)=0$, from which it follows that $\EE \|\Omega \xi(k)\|^2 \leq C \rho_{ess}(\bf{L}^T)\EE \|\Omega \xi(0)\|^2$ where $\rho_{ess}(\bf{L}^T)$ denotes the essential spectral radius of $\bf{L}^T$.
\end{proof}

\quad

Thanks to a result obtained in \cite{FF-SZ:08a}, the quantities $\mathbbm{E}\left[\|\xi -\alpha \1\|^2\right]$ and $\mathbbm{E}\left[\|\Omega \xi\|^2\right]$ have the same exponential convergence rate to zero, or, in other words, for any initial condition $\xi(0)$,
$$\limsup_{k\rightarrow\infty}\mathbbm{E}\left[\|\xi(k) -\alpha \1\|^2\right]^{1/k}=\limsup_{k\rightarrow\infty}\mathbbm{E}\left[\|\Omega \xi(k)\|^2\right]^{1/k}.$$ 
For this reason, in what follows we study the right-hand expression, which turns
out to be simpler to analyze. In order to have a single figure not dependent on the
initial condition, we focus on this worst case exponential rate of convergence
$$
R=\sup_{\xi(0)} \limsup_{k\rightarrow\infty}\mathbbm{E}\left[\|\Omega \xi(k)\|^2\right]^{1/k}
$$
It has been proved in Proposition $4.4$ of \cite{FF-SZ:08a} that
\begin{equation}\label{eq:rateBounds}
esr(\bar{Q})^2 \leq R \leq sr(\mathbb{E}(Q_i^T\Omega Q_i)).
\end{equation}
where $\bar{Q}$ denotes $\mathbbm{E}[Q_i]=\sum_{i=1}^N \beta_i Q_i$ and where $sr(\mathbb{E}(Q_i^T\Omega Q_i))$ denotes the spectral radius of semidefinite positive matrix $\mathbb{E}(Q_i^T\Omega Q_i)$, i.e., its largest eigenvalue. 

Unfortunately, it turns out from a numerical inspection over significant families of graphs, like Cayley graphs (see \cite{RC-FF-AS-SZ:08}), random geometric graphs, that the upper bound $sr(\mathbb{E}(Q_i^T\Omega Q_i))$ is greater than $1$, that is, it is not informative for our analysis. 

However we have run a number of MonteCarlo simulations randomized over graphs of different topology and size and over different initial conditions, and it always resulted that $\limsup_{k\rightarrow\infty}\mathbbm{E}\left[\|\Omega \xi(k)\|^2\right]^{1/k} \leq esr(\bar{Q})$. Based on this experimental evidence we formulate the following conjecture.

\begin{conjecture}
\label{conj:rateBound}
The quantity $esr(\bar{Q})$ is an upper bound for the exponential convergence rate $R$, i.e.,
$$
R \leq esr(\bar{Q}) 
$$
\end{conjecture}
The above conjecture and the fact that $esr(\bar{Q})^2 \leq R$ motivates to study $esr(\bar{Q})$. 

\section{ Convergence Rate Analysis of a-CL algorithm for regular graphs}
\label{sec:RateAnalysis}

In this section we assume that the measurements graph $\mathcal{G}=(V, \mathcal{E})$ is a connected undirected regular graph\footnote{An undirected graph is said to be regular if all the nodes have the same number of neighbors} such that, for $i\in \until{N}$, $|\mathcal{N}_i|=\nu$. In other words $\nu$ represents the degree of the graph. Moreover we assume the following properties.
\begin{assumption}\label{ass:1}
We have that
\begin{enumerate}
\item the error measurements covariances are all identical, i.e., $R_{ij}=R$ for all $(i,j)\in \mathcal{E}$;
\item  $\epsilon= R/(2(\nu+1))$;
\item the probabilities $\left\{\beta_1,\ldots,\beta_N\right\}$ are uniform, i.e., $\beta_1=\ldots=\beta_N=1/N$.
\end{enumerate} 
\end{assumption}
Observe that, from properties (i) and (ii) of the above assumption, it turns out that the matrix $P=I-\epsilon A^TR^{-1}A$, associated to the s-CL algorithm, is a symmetric matrix such that $P_{ij}=1/(\nu+1)$ for all $j\in \mathcal{N}_i \cup \left\{i\right\}$.\\

Let $\lambda_1(P)=1>\lambda_2(P)\geq\ldots\geq\lambda_N(P)$ the eigenvalues of $P$. Note that $esr(P)=\max \left\{|\lambda_2(p)|, |\lambda_N(P)|\right\}$. The goal is to show how the $2N$ eigenvalues of $\bar{Q}$ are related to those of $P$. 

\begin{lemma}
\label{lemma:eigenMap}
Consider the a-CL algorithm running over  an undirected regular graph $\mathcal{G}=(V, \mathcal{E})$ such that, for $i\in \until{N}$, $|\mathcal{N}_i|=\nu$. Assume Assumption \ref{ass:1} holds true.
Then the $2N$ eigenvalues of $\bar{Q}$ are the solutions of the following $N$ second-order equations
\begin{equation}
\label{eq:secondOrder}
f(s,\lambda_i,N,\nu) = s^2 + (a+b) s+ (ab + c)
\end{equation}
where
\begin{equation*}
\begin{split}
a &= - \left[ \frac{N-\nu}{N} + \frac{\lambda_i}{N}+\frac{\nu - 1}{N(\nu + 1)} \right] \\
b &= - \frac{N-1}{N}  \\
c &= -\frac{\nu-1}{N^2}(\lambda_i-\frac{1}{\nu+1})\\
\end{split}
\end{equation*}
\end{lemma} 
The proof is provided in appendix \ref{app:proofeigenMap}. Now let $s_1^{(i)}$ and $s_2^{(i)}$ denote the two solutions of $f(s,\lambda_i,N,\nu)$. It easy to see that $s_1^{(1)} = 1$ and $s_2^{(1)} = 1- \frac{\nu^2+1}{N(\nu+1)}$. The following result restricts the search of $esr(\bar{Q})$ among the values $|s_1^{(2)}|$ and $|s_2^{(2)}|$ and  $1- \frac{\nu^2+1}{N(\nu+1)}$.

\begin{theorem}\label{thm:LowerBoundRate}.  
Consider the a-CL algorithm running on  an undirected regular graph $\mathcal{G}=(V, \mathcal{E})$ such that, for $i\in \until{N}$, $|\mathcal{N}_i|=\nu$. Assume Assumption \ref{ass:1} holds true. Moreover let\\
$
\gamma^* = \frac{\nu-1+N(\nu+1)-\sqrt{N^2(\nu+1)^2 -2N(\nu^3+\nu+2)+(\nu-1)^2+(\nu^2+1)^2}}{\nu+1}
$
\\
then
\begin{enumerate}
\setlength{\itemsep}{5pt}
\item if $1-esr(P) \leq \gamma^* \Longrightarrow esr(\bar{Q})= \max (|s_1^{(2)}|, |s_2^{(2)}|)$;
\item if $1-esr(P) > \gamma^* \Longrightarrow esr(\bar{Q}) = s_2^{(1)} =  1- \frac{\nu^2+1}{N(\nu+1)}$.
\end{enumerate}
\end{theorem}
The proof is reported in appendix \ref{app:proofesrQ}.

\subsection{Asymptotic Behavior}
Consider now a sequence of connected undirected regular graphs $\mathcal{G}_N$ of increasing size $N$, and fixed degree $\nu$. Assume Assumption \ref{ass:1} holds true for any $\mathcal{G}_N$. Then to any $\mathcal{G}_N$ we can associate a stochastic matrix $P_N$ describing a s-CL algorithm, such that $(P_N)_{ij}=1/(\nu+1)$ for all $j\in \mathcal{N}_i \cup \left\{i\right\}$. 
In the following we analyze three different scenarios concerning three different spectral properties.\\

\subsubsection{Cayley graphs} The following assumption holds true for family of graphs like the Cayley graphs and $d$-dimensional  toruses. Moreover, it well describes the spectral behavior of Random Geometric graphs which is a family of graphs frequently used in practice.  

\begin{assumption}\label{ass:2}
Consider the sequence of matrices $P_N$ associated to the sequence of graphs $\mathcal{G}_N$ above described. Then we have 
\begin{equation}\label{eq:esrP}
esr\left(P_N\right) = 1-\varepsilon(N) + o(\varepsilon(N))
\end{equation}
\end{assumption}
where $\varepsilon: \N \to \R$ is a positive function such that $\varepsilon(N)\rightarrow0$ as $N\to \infty$. Let the matrix $\bar{Q}_N$ represent the average matrix associated to the a-CL algorithm running over $\mathcal{G}_N$. The following result characterizes the asymptotic behavior of $esr\left(\bar{Q}_N\right)$, with the respect to $esr\left(P_N\right)$.

\begin{theorem}\label{thm:proofAsymptotic}
Consider the sequence of graphs $\mathcal{G}_N$ above described. Consider the a-CL algorithm running over $\mathcal{G}_N$. Assume Assumption \ref{ass:1} and Assumption \ref{ass:2} hold true. Then
$$
esr\left(\bar{Q}_N\right) = 1 - \frac{\nu(\nu+1)}{N(\nu^2+1)}\varepsilon(N) + o\left(\frac{\varepsilon(N)}{N}\right)
$$
\end{theorem}
The proof can be found in appendix \ref{app:proofAsymptotic}.\\

\subsubsection{Ramanujan graphs} A similar analysis can be provided also for other relevant families of regular graphs like the Ramanujan graphs \cite{MRM:03}. Let us recall the asymptotic lower bound proved by Alon and Boppana for doubly stochastic matrices built over $\nu$-regular bidirected graphs. Specifically, if 
$A$ denotes the adjacency matrix of a $\nu$-regular bidirected graph, let $P$ be the doubly stochastic matrix defined as $P=\nu^{-1}A$, then
$$
\liminf_{N \to \infty} \rho_{\text{ess}}(P)\geq \frac{2 \sqrt{\nu-1}}{\nu}
$$
where the lim inf is taken along the family of all $\nu$-regular bidirected graphs having $N$ vertices. 
Ramanujan graphs are those $\nu$-regular bidirected graphs which achieves the previous bound, i.e., such that $\rho_{\text{ess}}(P)= \frac{2 \sqrt{\nu-1}}{\nu}$. Hence through the Ramanujan graphs it is possible to keep the essential spectral radius bounded away from $1$, while keeping the degree fixed. 

\begin{theorem}\label{thm:proofRamanujan}
Consider a sequence of Ramanujan graphs $\mathcal{G}_N$ of increasing size $N$ and fixed degree $\nu$. Consider the a-CL algorithm running over $\mathcal{G}_N$. Assume assumption \ref{ass:1} holds true. Then, it holds 
$$
\rho_{ess}(\bar{Q}_N)=1-\alpha(\nu)/N
$$
where, respectively\\ 
$\alpha(\nu) = \left[\frac{(\nu^2+1)+(1-\epsilon)(\nu+1)}{2(\nu+1)}\left(1-\sqrt{1-\frac{4\nu(1-\epsilon)(\nu+1)^2}{\left((\nu^2+1)+(1-\epsilon)(\nu+1)\right)^2}} \right)\right]$
\\
$
\epsilon = \frac{2}{\nu}\sqrt{\nu-1}
$
\end{theorem}

The proof can be found in Appendix \ref{app:proofRamanujan}.\\

\subsubsection{Complete graphs} Concerning the family of complete graphs we have that $\rho_{\text{ess}}(P)=0$ and $\nu=N-1$ so, we can state the following 
\begin{corollary}\label{thm:proofComplete}
Consider a sequence of Complete graphs $\mathcal{G}_N$ of increasing size $N$. Consider the a-CL algorithm running over $\mathcal{G}_N$. Assume assumption \ref{ass:1} holds true. Then, it holds 
$$
\rho_{ess}(\bar{Q}_N)=\frac{2(N-1)}{N}
$$ 
\end{corollary}
The proof can be found in Appendix \ref{app:proofComplete}.

\section{Numerical Results}
In this section we firstly provide some simulations supporting the exponential rate of convergence of the different algorithms implemented. Secondly some simulations to support conjecture \ref{conj:rateBound}.

\begin{example} 
In this example we consider a random geometric graph generated by choosing $N = 100$ points randomly placed in the interval $[0,10]$, and pairs of nodes that are within a range of $3$ took measurements of each others' relative positions, namely, both measurements $z_{ij}$ and $z_{ji}$ are available provided that $|x_i-x_j| \leq 3$.   
Every measurement was corrupted by Gaussian noise with unitary covariance. In this example we show the performance of a-CL with both the uniform persistent, assuming that each node transmitted at least once within each time interval $[k, k+\tau)$ where $\tau=200$, and the random persistent communicating network.

In Figure \ref{Fig:1} we plotted the behavior of the error
$$
J(k)= \log \left(\|A(\hat{x}(k)-x^*)\|\right)
$$
The plot reported is the result of the average over $1000$
Monte Carlo runs, randomized with respect to both the measurement graph\footnote{In perfoming our average we kept only the random geometric graphs which resulted to be connected.}
and the initial conditions. Observe that the trajectories converge to $0$ exponentially.

\begin{figure}[h]
\begin{center}
\begin{tabular}{c}
\includegraphics[width=1\columnwidth]{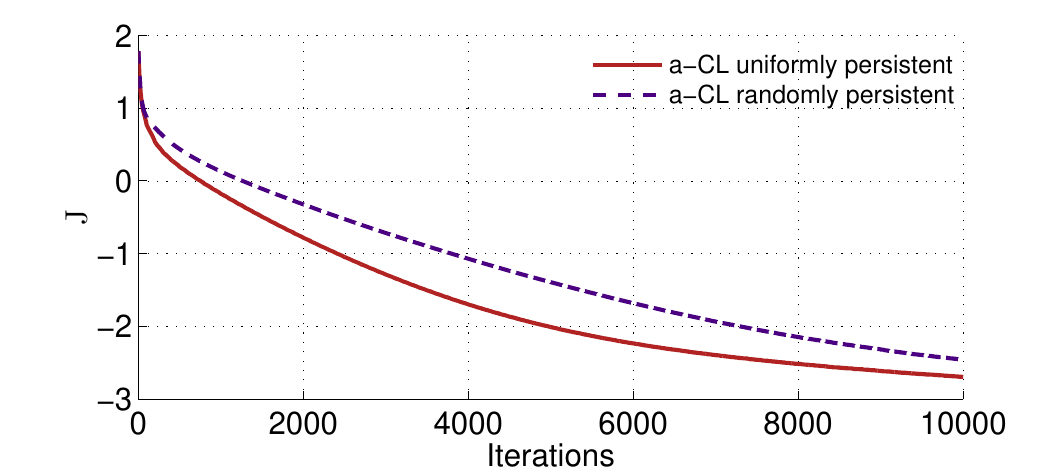}
\end{tabular}
\caption{\label{Fig:1} Behavior of $J$ for a uniformly persistent (red solid line) and a randomly persistent (blue dashed line) communicating network.}
\end{center}
\end{figure}

\end{example} 

%

\begin{example}
In this example we provide a numerical comparison with the randomized algorithm proposed in \cite{Bolognani:10}, later on called a-GL, which, for the sake of the completeness, we shortly review next. Similarly to the a-CL algorithm, also  during the $k$-th iteration in the randomized algorithm in \cite{Bolognani:10}, there is only one node, say node $i$, which transmits its variable $\hat{x}_i$ to all its neighbors. For $j \in \mathcal{N}_i$, node $j$, based on the information received from node $i$, performs the following update
\begin{eqnarray}\label{eq:Bolognani}
\hat{x}_j(k+1) &=&1/2\left(\hat{x}_j(k)+\hat{x}_i(k)+ 1/2(z_{ji}-z_{ij})\right)\\
&=&\hat{x}_j(k)+ 1/2\left(\hat{x}_j(k)-\hat{x}_i(k)+ 1/2(z_{ji}-z_{ij})\right) \nonumber
\end{eqnarray}
while for $h \notin \mathcal{N}_i$, $\hat{x}_h(k+1)=\hat{x}_h(k)$. 
In this example we consider random geometric measurement graphs $\mathcal{G}$ built as in the previous example. For both strategies the transmitting node at each iteration is chosen with a uniform probability. 
In Figure \ref{Fig:3} we plotted the behavior of $J$ for both the a-CL algorithm (blue dashed line) and the a-GL proposed in \cite{Bolognani:10} (red solid line).

As stressed in \cite{Bolognani:10} the trajectory $\hat{x}(k)$ generated by the updating rule in \eqref{eq:Bolognani} does not converge to $x^*$ but oscillates within a neighborhood of $x^*$. This fact appears to be evident also in Figure \ref{Fig:3}. However, of note is the fact that the strategy in \cite{Bolognani:10} seems to have a faster transient than the randomized a-CL algorithm. 
\begin{figure}
\begin{center}
\begin{tabular}{c}
\includegraphics[width=1\columnwidth]{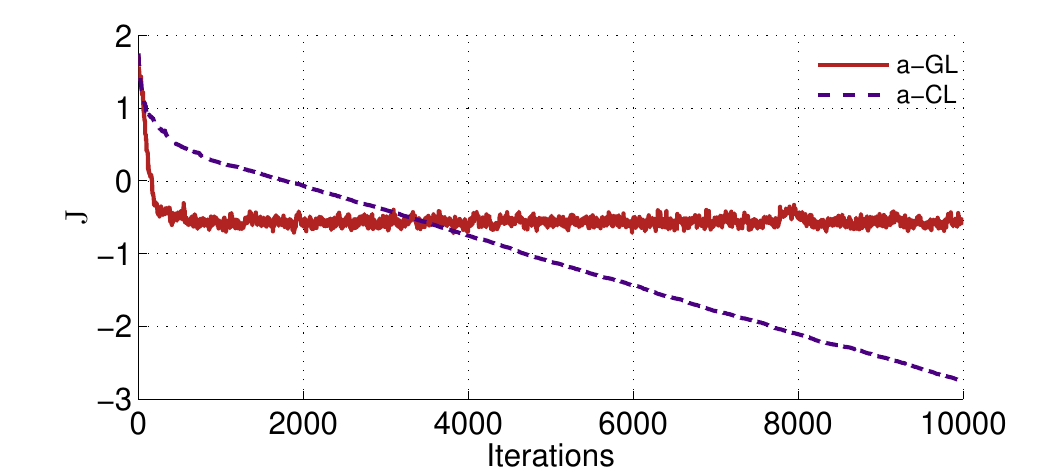}
\end{tabular}
\caption{\label{Fig:3} Comparison between the a-CL algorithm and the randomized algorithm (a-GL) proposed in \cite{Bolognani:10}.}
\end{center}
\end{figure}

\begin{remark}
The authors in \cite{Freris:12} proposed an algorithm (denoted as Randomized Kaczmarz Smoothing) similar to the one introduced in \cite{Bolognani:10}. However also the Randomized Kaczmarz Smoothing algorithm does not converge to the optimal estimate $x^*$ but oscillates within a neighborhood of it.
\end{remark}

\begin{remark}
The convergence properties of the randomized algorithm in \eqref{eq:Bolognani} might be improved by using the following modified updating rule
\begin{equation*}
\hat{x}_j(k+1) =\hat{x}_j(k)+ \gamma(k) \left(\hat{x}_j(k)-\hat{x}_i(k)+ 1/2(z_{ji}-z_{ij})\right) \nonumber
\end{equation*}
where $\left\{\gamma(k)\right\}$ is a sequence of receding step-sizes with $\gamma(k) \in (0,1)$ and $\lim_{k \to \infty} \gamma(k)=0$. A similar idea has been proposed in \cite{Freris:12} introducing the Randomized Kaczmarz algorithm with Under-Relaxation. However, in this case, since $\lim_{k \to \infty} \gamma(k)=0$ and, in turn, the applied control action becomes smaller and smaller, the speed of convergence is not exponential but just sub-linear. 
\end{remark}
\end{example}

\begin{example}
In this example we consider a set of Cayley graphs $\mathcal{G}$ of increasing size $N$. What we show is a comparison between the empirical rate of convergence of the algorithm, its lower bound, represented by $esr(\bar{Q})^2$, and the $esr(\bar{Q})$. As a matter of fact, figure \ref{fig:CayleyRate} shows that the rate is upper bounded by $esr(\bar{Q})$.
 
\begin{figure}[h]
\centering
\includegraphics[width=\columnwidth]{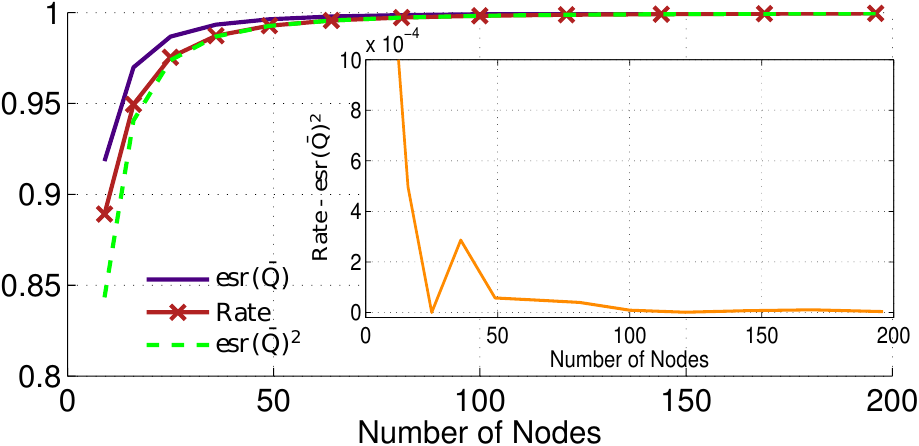}
\caption{Trend of the rate of convergence, of the $esr(\bar{Q})^2$ and of the $esr(\bar{Q})$ for Cayley graphs of increasing size $N$.}
\label{fig:CayleyRate}
\end{figure}

\end{example}

\begin{example}
In this final example we consider the family of random geometric graphs generated by choosing $N$ nodes in the interval [0,1] and connecting those whose relative euclidean distance is less than $0.3$. This meaning that node $i$ and $j$ perform the measurement $z_{ij}$ and $z_{ji}$ if $|x_i-x_j|\leq 0.3$. Figure \ref{fig:GeometricRate} shows that, for random geometric graphs as well as for Cayley, the rate is upper bounded by $esr{\bar{Q}}$.

\begin{figure}[h]
\centering
\includegraphics[width=\columnwidth]{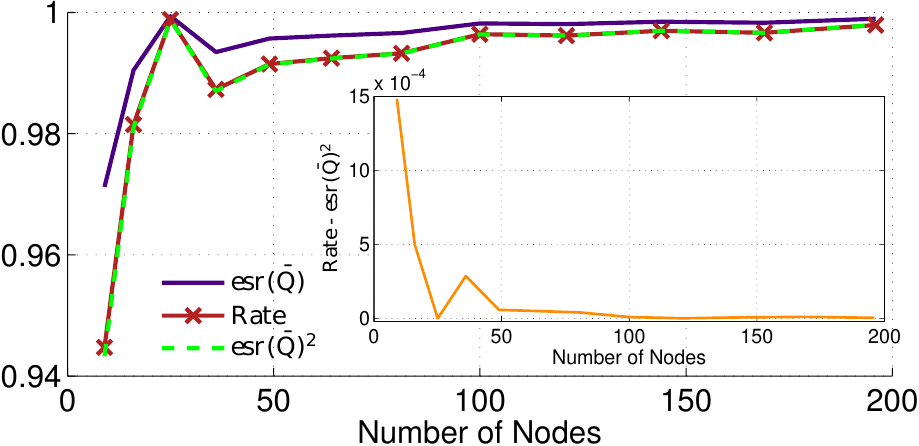}
\caption{Trend of the rate of convergence, of the $esr(\bar{Q})^2$ and of the $esr(\bar{Q})$ for random geometric graphs of increasing size $N$.}
\label{fig:GeometricRate}
\end{figure}

\end{example}

\begin{remark}
Note that, in both example presented, and in all the several simulations performed, the lower bound, $esr(\bar{Q})^2$, better approximates the behavior of the empirical rate of convergence.
\end{remark}


\section{Conclusion}
In this paper we considered the problem of optimal estimating the position of each agent in a network from relative noisy vectorial distances with its neighbors. We formulated the problem as a classical weighted least-squares problem and we reviewed its optimal centralized solution. The contribution we proposed is twofold.  We first  introduced a synchronous consensus-based algorithm that is guaranteed to have exponentially convergence rate to the optimal solution. We then extended this algorithm to more realistic asynchronous implementations  via the use of local memory variables. We showed that exponential convergence is still guaranteed under both uniform persistent communication and randomized persistent communication protocols.

Moreover we studied the rate of convergence of the proposed randomized broadcast asynchronous algorithm. We provided closed form expressions for the rate in terms of convergence in expectation for regular graphs.  We also conjectured based on extensive simulations that this provides also lower and upper bounds for the rate of convergence in mean square. We finally showed that asymptotically in the number of nodes $N$, this rate of convergence scales as standard memoryless consensus algorithms, thus implying a fast convergence rate. 

Future work is directed in proving the conjecture that the rate of convergence in expectation can provide upper and lower bounds for the rate of convergence in expectation, and in providing extensive comparisons with the most popular alternative algorithms available in the literature to solve the same problem.

\appendix

\subsection{Proof of Lemma \ref{lemma:eigenMap}}
\label{app:proofeigenMap}

Recalling equation (\ref{eq:Qi}), it is possible to write the average matrix $\bar{Q} = \mathbbm{E}[Q_i]$ as 
\[
\bar{Q} = \left[\begin{array}{c c}\bar{Q}_{11} & \bar{Q}_{12} \\ \bar{Q}_{21} & \bar{Q}_{22}\end{array}\right]
\]
where
\begin{eqnarray}
\bar{Q}_{11} &=& \left[\begin{array}{ccc}\frac{N-|\mathcal{N}_{1}|}{N} &  &  \\ & \ddots & \\ &  & \frac{N-|\mathcal{N}_{N}|}{N}\end{array}\right]I  +  \frac{1}{N}P \nonumber  \\ 
&+& \left[\begin{array}{ccc}\frac{|\mathcal{N}_{1}|-1}{N} &  &  \\ & \ddots &  \\ &  &\frac{|\mathcal{N}_{N}|-1}{N} \end{array}\right]diag(P)\nonumber\\
 \bar{Q}_{12} &=& \left[\begin{array}{ccc}\frac{|\mathcal{N}_{1}|-1}{N} &  &  \\ & \ddots &  \\ &  & \frac{|\mathcal{N}_{N}|-1}{N}\end{array}\right]\left(P-diag(P)\right)\nonumber\\ 
\bar{Q}_{21} &=& \frac{1}{N}I \nonumber\\
\bar{Q}_{22} &=& \frac{N-1}{N}I\nonumber
\end{eqnarray}

Thanks to assumption (i) we have that $$\mathcal{N}_{i}=\nu \qquad \forall i \in \until{N}$$
therefore we can write $\bar{Q}$ in a more compact form, where

\begin{equation*}
\begin{split}
&\bar{Q}_{11} =   \frac{N-\nu}{N}I + \frac{1}{N}P + \frac{\nu-1}{N(\nu+1)}I \\
&\bar{Q}_{12} = \frac{\nu-1}{N} (P-\frac{1}{\nu+1}I)\\
&\bar{Q}_{21} = \frac{1}{N}I \\
&\bar{Q}_{22} = \frac{N-1}{N}I
\end{split} 
\end{equation*}

Since $P$ and $I$ are simultaneously diagonalizable, with an appropriate permutation is possible to obtain a block diagonal matrix where every block, of dimension two, is equal to
\[
\begin{bmatrix} 
\frac{N-\nu}{N} + \frac{\lambda_i}{N} + \frac{\nu-1}{N(\nu+1)} &
\qquad\frac{\nu-1}{N} (\lambda_i-\frac{1}{\nu+1})\\
\frac{1}{N} &
\frac{N-1}{N}
\end{bmatrix}
\]
where $\lambda_i$ represents the generic eigenvalue of $P$. Thanks to \cite{Silvester:00} the matrix entire determinant can be obtain as the multilpcation of the single block determinant whose associated characteristic polynomial is equal to:
\begin{eqnarray}\label{eq:polyCar}
f(s,\lambda_i,N,\nu) &=& s^2 - \Bigl(\frac{N-\nu}{N} + \frac{\lambda_i}{N}+\frac{\nu - 1}{N(\nu + 1)}+\nonumber\\
								&+&\frac{N-1}{N}\Bigr)s+\nonumber\\
								&+& \frac{N-1}{N}\left(\frac{N-\nu}{N} + \frac{\lambda_i}{N}+\frac{\nu - 1}{N(\nu + 1)}\right)-\nonumber\\ 
								&-&\frac{\nu-1}{N^2}\left(\lambda_i-\frac{1}{\nu+1}\right)
\end{eqnarray}
Thanks to (\ref{eq:polyCar}) it can be seen how every eigenvalue of $P$ is mapped in two eigenvalues of $\bar{Q}$.

\subsection{Proof of Theorem \ref{thm:LowerBoundRate} } 
\label{app:proofesrQ}

The proof of the theorem arises from a collection of preliminary results that we are now going to show. It is known that, being $P$ a stochastic matrix, all its eigenvalues can be written like $\lambda_i = 1 - \gamma_i$, where $\gamma_i \in [0,2]$. Now rewriting the second order equation (\ref{eq:secondOrder}) in two parts, one depending from $\gamma_i$ and one not depending from it, we obtain:
\begin{equation}\label{eq:polyCarGamma}
f(s,\lambda_i,N,\nu) = d(s,N,\nu) + \gamma_in(s,N,\nu)
\end{equation}
Note that, this can be considered as a "quantized" root locus, where the gain is given by the N values of $\gamma_i$. For sake of simplicity we will study its continuous time version in order to understand which are the possible eigenvalues of $\bar{Q}$.
\begin{remark}
It is known that the eigenvalues are discrete but, to the end of exploiting the root locus technique, we  will consider their "continuous" version. We want to highlight that even if we have infinite values of $\gamma_i$ in the root locus, we are interested in just 2N points on it.
\end{remark}

As shown in equation (\ref{eq:polyCarGamma}) the two polynomials involved in the locus are 
\begin{equation*}
\begin{split}
d(s,N,\nu) &= s^2 - \frac{2N(\nu+1)-(\nu^2+1)}{N(\nu+1)}s + 1- \frac{\nu^2+1}{N(\nu+1)}\\
n(s,N,\nu) &= \frac{s}{N} + \frac{\nu-N}{N^2}
\end{split}
\end{equation*}
whose roots are $d_1 = 1$, $d_2 = 1- \frac{\nu^2+1}{N(\nu+1)}$ and $n_0 = 1-\frac{\nu}{N}$. We recall from standard results that $d_1$ and $d_2$ are the starting point of the root locus\footnote{note that they coincide with the solution of equation (\ref{eq:secondOrder}), previously called $s_1^{(1)}$ and $s_2^{(1)}$}, and $n_0$ its ending point, respectively for $\gamma = 0$ and $\gamma = \infty$. Additionally, the double points of the locus are 
$$
dp_{+,-} = 1- \frac{1}{N}\left(\nu \pm \sqrt{\frac{\nu(\nu-1)}{\nu+1}}\right)
$$


\begin{figure}%
\centering
\subfloat[][]{\label{fig:rootLocus}\includegraphics[width=0.40\columnwidth]{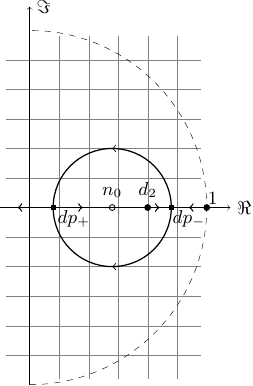}}%
\subfloat[][]{\label{fig:vectors}\includegraphics[width=0.40\columnwidth]{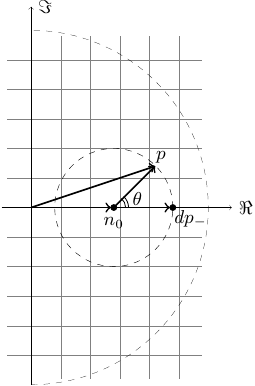}}%

\caption{(a) Root Locus of $d(s) + \epsilon_i n(s) = 0$. (b) The figure shows that the length of the vector between the points $0$ and $dp_{-}$ is   longer than the vector formed by the points $0$, $n_0$ and $p$}%
\label{fig:cont}%
\end{figure}

It is easy to see that
$$
0 < n_0 \leq d_2 \leq d_1 
$$
indeed, the first inequality is trivial being $\nu < N$, while the second leads to
$$
\frac{N-\nu}{N} - 1 + \frac{\nu^2+1}{N(\nu+1)} \leq 0 \Rightarrow \nu \geq 1 
$$
Moreover is straightforward to see that $n_0$ is the mid point of $dp_+$ and $dp_-$ being
\begin{equation*}
\begin{split}
n_0 =& 1-\frac{\nu}{N} = \frac{dp_- + dp_+}{2} =\\
=& \frac{1}{2}\left( 1-\frac{1}{N}\left(\nu - \sqrt{\frac{\nu(\nu-1)}{\nu+1}}\right)\right)+  \\
&+ \frac{1}{2}\left(1-\frac{1}{N}\left(\nu + \sqrt{\frac{\nu(\nu-1)}{\nu+1}}\right)\right)
\end{split}
\end{equation*}

It follows that the double points are centred on $n_0$ and that $dp_-$ lays between $d_1$ and $d_2$.\\

Let us firstly show the following two corollaries.  

\begin{corollary}\label{cor:rootLocusCircle}
Given a root locus function of the form\footnote{This kind of equation results from the root locus of a transfer function with one zero and two poles.} $(s-p_1)(s-p_2)+K(s-z) = 0, \quad K \in [0,+\infty)$,  the complex roots lay on a circle centred on the $z$.
\end{corollary}
\begin{proof}
The roots of $(s-p_1)(s-p_2)+K(s-z) = 0$ are:
$$
s_{1,2} = \frac{(p_1+p_2-K) \pm \sqrt{(p_1+p_2-K)^2-4(p_1p_2-Kz)}}{2}
$$
whose real and imaginary part, considering the complex case that holds when $(p_1+p_2-K)^2-4(p_1p_2-Kz) < 0$,
are:
$$
\Re(s_{1,2}) = \frac{p_1+p_2-K}{2}
$$
$$
\Im(s_{1,2}) = \pm \frac{\sqrt{(p_1+p_2-K)^2-4(p_1p_2-Kz)} }{ 2 }
$$

The square of the distance of these points from $z$ can be computed as
$$
(\Re(s_{1,2})-z)^2 + (\Im(s_{1,2}) - 0)^2 = z^2 - zp_1 - zp_2 + p_1p_2
$$
that is independent on $K$ and so is a constant.
\end{proof}

Thank to corollary (\ref{cor:rootLocusCircle}), is easy to get that the absolute value of the complex roots is descending when $\gamma_i$ increase.


In particular every pair of complex roots has module smaller that $dp_-$, indeed, being $|dp_{-}|-|n_0|$ the radius of the circonference centred in $n_0$, see figure \ref{fig:vectors}, it follows that
\begin{eqnarray}
   |p|^2 &=& (|n_0| + (|dp_{-}| - |n_0|)\cos\theta)^2 + \nonumber\\
   			&+& ((|dp_{-}| - |n_0|)\sin\theta)^2 \nonumber \\
   &=& |dp_{-}|^2 + 2|n_0|(1-\cos\theta)(|n_0|-|dp_{-}|)\leq |dp_{-}|^2 \nonumber
\end{eqnarray}

since $(1-\cos\theta) \geq 0$ and $(|n_0|-|dp_{-}|) \leq 0$. Where $\theta$ is the angle of the vector $p-n_0$ wrt the real axis.\\

\begin{corollary}\label{cor:Re>0}
Consider polynomial in equation (\ref{eq:polyCarGamma}). For $N>1$ and $\nu \in \until{N-1}$ then $\forall$ $\gamma_i \in [0,2]$, the real parts of $s_{1,2}^{(i)}$, solution of $f(s,\lambda_i,N,\nu)$, are non negative.    
\end{corollary}
\begin{proof}
The generic solution of the characteristic polynomial (\ref{eq:polyCarGamma}) are equal to
$
s_{1,2}^{(i)} = \frac{1}{2}\left(\frac{2N(\nu+1)-(\nu^2+1)}{N(\nu+1)} - \frac{\gamma_i}{N}\right) \pm \frac{1}{2}\sqrt{\left(\frac{2N(\nu+1)-(\nu^2+1)}{N(\nu+1)} - \frac{\gamma_i}{N}\right)^2 - 4\left(1-\frac{\nu^2+1}{N(\nu+1)} + \gamma_i\frac{\nu-N}{N^2}\right)} 
$
whose discriminant could be $\Delta \lesseqgtr 0 $. Let us firstly consider $\Delta \geq 0$. In this case the value of $\gamma_i$ such that $f(s,1-\gamma_i,N,\nu)$ has one root equal to zero is 
$$
\gamma_i = \frac{N(\nu^2+1-N\nu-N)}{(\nu+1)(\nu-N)}
$$
Imposing $\gamma_i \geq 2$, after some algebra, one gets
$$
(N-2)\nu^2 + (2N-2-N^2)\nu + (3N-N^2) \leq 0
$$
which is always satisfied for $N>1$ and $\nu \in \until{N-1}$.\\

Consider now, the second case with $\Delta<0$. It is easy to see that the two pure imaginary roots of $f(s,1-\gamma_i,N,\nu)$ are reached with a value of $\gamma_i\geq 2$. Indeed, imposing $\Re(s_{1,2}^{(i)}) = 0$ ones get
$$
\gamma_i = 2N - \frac{\nu^2+1}{\nu+1}
$$
which, $\forall N\geq 2$, is a descending function of $\nu$ and assumes its minimum value in $\nu=1$ which is
$$
\gamma_i = 2N - 1 \geq 2
$$

\end{proof}

Thanks to corollary \ref{cor:rootLocusCircle} and \ref{cor:Re>0} we know that, $\forall \lambda_i(P)$, all the eigenvalues of $\bar{Q}$ have positive real part. Moreover analyzing the root locus form it is easy to understand that three main situations can occur, that are:
\begin{enumerate}
\item the value of $\gamma_i$ corresponding to $esr(P)$ is sufficiently small to individuate, on the root locus, two eigenvalue of $\bar{Q}$ that lie on the real axis between $1$ and $d_2$. In this case the greater correspond to the $esr(\bar{Q})$;
\item the value of $\gamma_i$ of the $esr(P)$ individuates, on the locus, two complex conjugate values that has module greater than $d_2$. In this case, again, the $esr(\bar{Q})$ correspond to them;
\item $\gamma_i$ individuates values that have module smaller than $d_2$. In this case is obvious that $d_2$ becomes the $esr(\bar{Q})$
\end{enumerate}

This is enough to conclude that, calling $\gamma^*$ the value that individuates two complex conjugate roots on the locus, say $\mu_{1,2}$, such that 
$$
|\mu_1| = |\mu_2| = |d_2|
$$
then, if $\gamma<\gamma^*$ we get that the $esr(P)$ maps the $esr(\bar{Q})$; otherwise the $esr(\bar{Q})$ is equal to $d_2$. The critical value of $\gamma^*$ is\\
\\
$
\
\gamma^* = \frac{\nu-1+N(\nu+1)-\sqrt{N^2(\nu+1)^2 -2N(\nu^3+\nu+2)+(\nu-1)^2+(\nu^2+1)^2}}{\nu+1}
$

This complete the proof.

\subsection{Proof of Theorem \ref{thm:proofAsymptotic}}
\label{app:proofAsymptotic}
\begin{proof}
Let us define the characteristic polinomial that identify the $2N$ eigenvalues of $\bar{Q}_N$ wrt those of $P_N$, denoted as $1-\gamma_i$, as 
$$
f(s,\lambda_i,N,\nu) = d(s,N,\nu) + \gamma_in(s,N,\nu) \triangleq g(s,\gamma_i,N,\nu)
$$
where $g$ explicity explicitly shown the dependance on $\gamma_i$ and where
\begin{equation*}
\begin{split}
d(s,N,\nu) &= s^2 - \frac{2N(\nu+1)-(\nu^2+1)}{N(\nu+1)}s + 1- \frac{\nu^2+1}{N(\nu+1)}\\
n(s,N,\nu) &= \frac{s}{N} + \frac{\nu-N}{N^2}
\end{split}
\end{equation*}
It can be easily shown that the $sr(P_N)$ maps the $sr(\bar{Q}_N)$, i.e choosing $\gamma_i=0$ one of the two solution of $g$ is $s=1$.\\
Thank to assumption \ref{ass:2} we have that 
$$
esr(P_N) = 1-\varepsilon(N) + o(\varepsilon(N)) 
$$ 
and since
\begin{equation*}
\left\lbrace
\begin{array}{c}
g(1,0,N,\nu) = 0\\
\left.\frac{\partial g}{\partial s}\right|_{(1,0,N,\nu)} \neq 0
\end{array}
\right.
\end{equation*}
it is possible to exploit the implicit function theorem that let us write, denoting with $s(\gamma_i)$ the eigenvalues of $\bar{Q}_N$ as function of the corresponding $\gamma_i$,
\begin{eqnarray}
s(\varepsilon(N)+o(\varepsilon(N))) &=& 1 - \left.\frac{\partial g}{\partial \gamma_i}\slash\frac{\partial g}{\partial s}\right|_{(1,0,N,\nu)}(\epsilon(N)+o(\epsilon(N)))\nonumber\\
&=& 1- \frac{\nu(\nu+1)}{N(\nu^2+1)}\epsilon(N) + o\left(\frac{\epsilon(N)}{N}\right)\nonumber
\end{eqnarray}

which means that $esr(\bar{Q})$ can be expressed as
\begin{equation}\label{eq:esrQbar}
esr(\bar{Q}) = 1 - \frac{\nu(\nu+1)}{N(\nu^2+1)}\epsilon(N) + o\left(\frac{\epsilon(N)}{N}\right)
\end{equation}

\end{proof}

\subsection{Proof of Theorem \ref{thm:proofRamanujan}}
\label{app:proofRamanujan}

\begin{proof}
Recalling that $$\lim_{N\to \infty}\gamma^*(\nu,N)=\frac{\nu^2+1}{(\nu+1)^2}$$ and being $\nu>2$,  can be easily proved that  $1- \rho_{\text{ess}}(\bar{P_N}) < \gamma^*$ so, thanks to Theorem \ref{thm:LowerBoundRate}, it follows that $\rho_{\text{ess}}(\bar{Q}_N)= \max (|s_1^{(2)}|, |s_2^{(2)}|)$. Moreover, after straightforward algebraic calculations it is easy to see that $|s_1^{(2)}| > |s_2^{(2)}|$, therefore 
\begin{equation}
\rho_{\text{ess}}(\bar{Q}_N) = s_1^{(2)} = 1 - \frac{\alpha(\nu)}{N}
\end{equation}
where, respectively\\
$
\alpha = \left[\frac{(\nu^2+1)+(1-\epsilon)(\nu+1)}{2(\nu+1)}\left(1-\sqrt{1-\frac{4\nu(1-\epsilon)(\nu+1)^2}{\left((\nu^2+1)+(1-\epsilon)(\nu+1)\right)^2}} \right)\right]
$\\
$
\epsilon = \frac{2}{\nu}\sqrt{\nu-1}
$
\end{proof}

\subsection{Proof of Corollary \ref{thm:proofComplete}}
\label{app:proofComplete}
\begin{proof}
Since all the eigenvalues of $P_N$ different from $1$ are equal to $0$ the essential spectral radius $\rho_{ess}(P_N)=0$ and the spectral gap $\gamma_2 = 1$. So it easy to see that

$$
1-\rho_{ess}(P_N) = 1 > \gamma^* = \frac{\nu^2+1}{(\nu+1)^2}
$$

this meaning that, for theorem (V.6)

$$
\rho_{ess}(\bar{Q}_N) = 1-\frac{\nu^2+1}{N(\nu+1)} = \frac{2(N-1)}{N^2}
$$
where the last equality follows from the fact that for a complete graph $\nu=N-1$.
\end{proof}

\bibliographystyle{IEEEtran}
\bibliography{bib/refs}

\end{document}